\documentclass{amsart}

\usepackage{comment}
\usepackage{amsmath}
\usepackage{amsfonts}
\usepackage{amssymb}
\usepackage{amsthm}
\usepackage{dsfont}
\usepackage{accents}
\usepackage[alphabetic, msc-links]{amsrefs}
\usepackage{mathtools}
\usepackage{xfrac}
\usepackage{hyperref}
\usepackage{float}
\usepackage{tikz-cd}
\usepackage{graphicx}
\usepackage[normalem]{ulem}
\usepackage{subfiles}

\newcommand{\Id}{\mathds{1}}

\newcommand{\RR}{\mathbb{R}}
\newcommand{\ZZ}{\mathbb{Z}}

\newcommand{\sph}{\mathbb{S}}

\newcommand{\DD}{\mathbb{D}}

\newcommand{\ecal}{\mathcal{E}}
\newcommand{\dcal}{\mathcal{D}}
\newcommand{\hcal}{\mathcal{H}}
\newcommand{\scal}{\mathcal{S}}
\newcommand{\ncal}{\mathcal{N}}
\newcommand{\pcal}{\mathcal{P}}
\newcommand{\ucal}{\mathcal{U}}

\newcommand{\Df}{\mathbf{D}}

\newcommand{\SO}{\operatorname{SO}}
\newcommand{\Eq}{\ecal q}
\newcommand{\Eqp}{\Eq_{+}}
\newcommand{\RP}{\mathbb{R}\mathbf{P}}

\newcommand{\Diffp}[1]{\operatorname{Diff}^+(#1)}
\newcommand{\Ham}[1]{\operatorname{Ham}(#1)}
\newcommand{\Hamc}[1]{\operatorname{Ham}_c(#1)}

\newcommand{\Symp}{\operatorname{Symp}}
\newcommand{\eval}{\operatorname{ev}}

\newcommand{\Area}{\operatorname{Area}}
\newcommand{\girth}{\operatorname{girth}}
\newcommand{\diam}{\operatorname{diam}}
\newcommand{\Stab}{\operatorname{Stab}}

\newtheorem*{remark*}{Remark}
\newtheorem{lemma}{Lemma}
\newtheorem{theorem}{Theorem}
\newtheorem*{theorem*}{Theorem}

\newtheorem*{definition*}{Definition}
\theoremstyle{remark}
\newtheorem*{acknowledgements*}{Acknowledgements}
\newtheorem*{structure}{Structure of the paper}

\title{On the Hofer Girth of the Sphere of Great Circles}
\author{Itamar Rosenfeld Rauch}
\address{Itamar Rosenfeld Rauch, Mathematics Department, Technion - Israel Institue of Technology, Haifa, 32000, Israel}
\email{itamar-ros@campus.technion.ac.il}
\thanks{The author was partially supported by the Azrieli Foundation.}

\begin{document}

\maketitle

\begin{abstract}
An oriented equator of $\sph^2$ is the image of an oriented embedding $\sph^1 \hookrightarrow \sph^2$ such that it divides $\sph^2$ into two equal area halves.
Following Chekanov, we define the Hofer distance between two oriented equators as the infimal Hofer norm of a Hamiltonian diffeomorphism taking one to another.
Consider $\Eqp$ the space of oriented equators.
We define the Hofer girth of an embedding $j:\sph^2 \hookrightarrow \Eqp$ as the infimum of the Hofer diameter of $j'(\sph^2)$, where $j'$ is homotopic to $j$.
There is a natural embedding $i_0:\sph^2\hookrightarrow\Eqp$,  sending a point on the sphere to the positively oriented great circle perpendicular to it.
In this paper we provide an upper bound on the Hofer girth of $i_0$.
\end{abstract}

\section{Introduction}
\subsection{Definitions and main result}
Let $\sph^2\subset \RR^3$ be the two dimensional unit sphere, equipped with a symplectic structure $\omega= \frac{1}{4\pi} \omega_{\mathrm{std}}$ such that $\Area(\sph^2)=1$. 
An \emph{(oriented) equator} $L$ of $\sph^2$ is the image of an (oriented) embedding $\sph^1 \to \sph^2$ such that it divides the sphere into two equal area halves. 
Denote by $\Eq$ the space of unoriented equators, and likewise by $\Eqp$ the space of oriented equators.
This paper is a part of an endeavor to better understand the geometry of the $\sph^2$ equators space.
While not a direct sequel to \cite{MR2597651}, and \cite{MR3395269}, this work is inspired by these articles in an attempt to investigate special Lagrangians on surfaces.

Consider a compact Riemannian manifold $(X, g)$. 
Recall that the \emph{systole} of $X$ is a metric invariant of $X$, defined as the least length of a noncontractible loop in $X$.
This sort of 1-dimensional question is classical in geometry.
Higher dimensional generalizations go back to Loewner, Berger, Pu, and others; see \cite{MR1427763} for a survey.
The work presented here may be considered as one way to generalize systoles to spheres in $\Eqp$.
Roughly speaking, we would like to find the least Hofer diameter of an embedded noncontractible sphere in $\Eqp$.

To define this Hofer diameter we need to recall the definition of Hofer distance; for more details see, e.g., \cite{MR3674984}. 
Let $(M,\omega)$ be a connected symplectic manifold without boundary. Denote by $\Hamc{M}$ the group of compactly supported Hamiltonian symplectomorphisms. The Hofer distance between $\phi_0,\phi_1\in \Hamc{M}$ is given by
\begin{equation*}
\rho(\phi_0, \phi_1) = \inf_{\phi_H = \phi_1 \circ \phi_0^{-1}} \int_{0}^{1} \max_{p\in M} H(t, p) - \min_{p \in M} H(t, p) dt,
\end{equation*}
where the infimum goes over all time dependent compactly supported smooth Hamiltonians $H:[0,1] \times M \to \RR$ that generate $\phi_1 \circ \phi_0^{-1}$.

In this paper we only consider $\sph^2$ as our symplectic manifold, and since it is compact $\Hamc{\sph^2} = \Ham{\sph^2}$. 
Consider the following distance function on $\Eqp$, introduced by Chekanov in \cite{MR1774099}.
\begin{definition*}
	Let $L_1,L_2\subset \sph^2$ be two oriented equators, define their distance by
	\begin{equation*}
	d_H(L_1,L_2) := \inf_{\substack{\phi\in \Ham{\sph^2} \\ \phi(L_1) = L_2 }} \Vert \phi \Vert.
	\end{equation*}
\end{definition*}
The Hofer distance induces the following definitions of diameter and girth.
\begin{definition*}
	Let $f:\sph^2\to \Eqp$. Put,
	\begin{equation}\label{eq:def-Hofer-diam-of-emb}
	\diam(f) := \sup_{x,y\in\sph^2} d_H(f(x), f(y)).
	\end{equation}
	Let $\alpha\in\pi_2(\Eqp)$. Define,
	\begin{equation}\label{eq:def-Hofer-diam-for-classes}
	\girth(\alpha) = \inf_{f\in\alpha} \diam(f).
	\end{equation}
\end{definition*}

Our main result is the following.
\begin{theorem}\label{thm:main-theorem}
	Let $i_0:\sph^2\to \Eqp$ denote the embedding sending each point on the sphere to the positively oriented great circle perpendicular to it. 
	Then, the following bound holds for the homotopy class $[i_0]\in\pi_2(\Eqp)$:
	\begin{equation*}
	\girth([i_0]) \leq 1/3.
	\end{equation*}
\end{theorem}
Note that the Hofer diameter of $i_0$ is $1/2$, as we explain in Subsection \ref{sec:Elem-observations}.
The proof of Theorem \ref{thm:main-theorem} consists in assembling the pieces from the following sections, for completeness we provide it here.
\begin{proof}
	
	Fix some sufficiently small $\delta,\epsilon>0$.
	We construct in Section \ref{sec:homotopy-construction} a homotopy depending on $\delta$, and $\epsilon$, between $i_0$ and some embedding $i_p$, which we call a ``pipe equator embedding".
	We calculate in Section \ref{sec:hofer-diam-calculation} that $1/3+\delta$ is an upper bound for the Hofer diameter of $i_p$.
	Since $\delta$ is arbitrary, Theorem \ref{thm:main-theorem} follows.
	
\end{proof}

Remark that there are other interesting attempts at understanding the geometry of $\Eqp$. 
We would like to note one such attempt by Y. Savelyev, somewhat close in spirit to the contents of this paper; See \cite{MR3356596} and Savelyev's preprint in \cite{savelyev2015global}\footnote{We thank Egor Shelukhin for referring us to this paper.}, and more advanced versions of this paper in the author's website.
Savelyev calculates a similar notion to the girth above in the space of loops of equators with respect to the positive Hofer length functional.

\subsection{Four elementary observations}\label{sec:Elem-observations}
First, the homotopy class of $i_0$ is non-trivial; see Section \ref{sec:algebraic-topology} for details. 
Hence, whether or not $\girth([i_0])$ is zero is a non-trivial question.
Second, we have the following naive upper bound on $\girth([i_0])$,
\begin{equation}\label{eq:naive-bound}
\diam(i_0) = \frac{1}{2}.
\end{equation}
We claim that the Hofer distance between any two great circles is bounded above by $1/2$.
Indeed, recall that $SO(3)$ is a subgroup of $\Ham{\sph^2}$; see, e.g. Section 1.4 in \cite{MR1826128}.
Thus the Hofer distance between any two great circles is bounded above by the Hofer norm of the rotation taking one to another.
Under the choice of symplectic form with $\Area(\sph^2)=1$, it holds that the Hofer norm of such rotations is bounded above by $1/2$.
To see this, let $R=R(\theta, u)$ be a counterclockwise rotation with angle $\theta$ about the axis prescribed by the unit vector $u$.
It is straightforward to check that the height function with respect to $u$, i.e., $F(x) = (x, u)$ defined for $x\in\sph^2$, where here $(\cdot, \cdot)$ stands for the standard Euclidean inner product, generates $R$ as a Hamiltonian diffeomorphism; see example 1.4.H in \cite{MR1826128}.
The time parameter of the flow generated by $F$ determines the angle $\theta$ for $R$.
Given an angle $\theta \in [0, \pi]$, it follows that
\begin{equation*}
\Vert R(\theta, u)\Vert \leq \int_{0}^{\theta} \max_{x\in \sph^2} F(x) - \min_{x\in\sph^2}F(x) dt = 2\theta.
\end{equation*}
Rescaling by $4\pi$ to take our choice of area form into account, yields
\begin{equation*}
\Vert R(\theta, u)\Vert \leq \frac{2\theta}{4\pi} \leq \frac{1}{2}.
\end{equation*}
This shows the upper bound in \eqref{eq:naive-bound}; and in particular,
\begin{equation*}
\girth([i_0]) \leq \frac{1}{2}.
\end{equation*}

The lower bound in \eqref{eq:naive-bound} follows from the well-known energy-capacity inequality, which takes the following simpler form for surfaces.
Given a closed symplectic surface $(M,\omega)$, for any connected subset $A\subset M$ it holds that $e(A)\footnotemark \geq \Area(A)$.
\footnotetext{Recall that the displacement energy $e(A)$ for a connected $A\subset (M,\omega)$ is defined to be $\inf\Vert \phi \Vert$, where the infimum goes over $\phi\in \Hamc{M}$ with $A\cap \phi(A) = \emptyset$. }
To apply it to our case, note that $i_0(x)$ bounds a disc $D\subset\sph^2$ of area $1/2$, and that $i_0(-x)$ is the same great circle as $i_0(x)$ with orientation reversed. 
Meaning, any Hamiltonian diffeomorphism $\phi$ taking $i_0(x)$ to $i_0(-x)$ would displace $D$, and so we find that
\begin{equation*}
\Vert \phi \Vert \geq e(D) \geq \Area(D) = \frac{1}{2}.
\end{equation*}
Hence, Equation \eqref{eq:naive-bound} holds.

The third elementary observation is that if we consider the space of unoriented equators $\Eq$, and the embedding $\iota:\RP^2 \to \Eq$ induced by $i_0$, we find that,
\begin{equation*}
\girth([\iota]) = \diam(\iota) = 1/4.
\end{equation*}
By this we mean that the Hofer diameter of $\iota$ is $1/4$, and maps homotopic to $\iota$ have Hofer diameter of at least $1/4$.
To define $\iota$ precisely, let $q\in \RP^2$; lift it to $\tilde{q}\in\sph^2$ via the double cover $\sph^2\to\RP^2$.
Then, put $\iota(q)$ to be the great circle $i_0(\tilde{q})$ without the orientation.
Due to the composition with the orientation forgetful map $\Eqp \to \Eq$, it holds that $\iota(q)$ does not depend on the lift $\tilde{q}$.

We claim now that $d:=\diam(\iota)=1/4$.
For ease of notation, consider the induced Hofer distance on $\RP^2$, defined by
\begin{equation*}
d_H(x,y) = d_H(\iota(x), \iota(y)), \quad \forall x,y\in \RP^2.
\end{equation*}
So, $d$ is the diameter of $\RP^2$ under this distance, i.e.,
\begin{equation*}
d =\diam (\RP^2) := \sup_{x,y\in\RP^2} d_H(x, y).
\end{equation*}
To prove the claim, observe that an upper bound for $d$ is obtained by a rotation by $\pi/2$, similar to the case of the second elementary observation above.
For the lower bound, note that for each point $q\in\RP^2$ there exists another point $p\in\RP^2$, such that
\begin{equation*}
d_H(q, p) \geq 1/4.
\end{equation*}
To see this, consider any non-contractible loop in $\RP^2$ based in $q$. 
It lifts to a curve $\gamma$ in $\sph^2$ connecting two antipodal points $x,-x$, where $q$ is the class of $\{x,-x\}$ under the covering map $\sph^2 \to \RP^2$.
By \eqref{eq:naive-bound}, it follows that $d_H(i_0(x), i_0(-x)) \geq 1/2$, whence the Hofer length of $\gamma$ is at least $1/2$.
Moreover, there exists a point $\tilde{p}\in\sph^2$ lying on $\gamma$, such that $d_H(i_0(\tilde{p}), i_0(x)) \geq 1/4$, and likewise $d_H(i_0(\tilde{p}), i_0(-x)) \geq 1/4$.
Projecting $\tilde{p}$ to $\RP^2$ yields a point $p\in\RP^2$ with
\begin{equation*}
d_H(q, p) \geq 1/4.
\end{equation*}

Let $H:\RP^2\times I \to \Eq$ be a homotopy of embeddings $\RP^2 \to \Eq$, where $I=[0,1]$ and $H_0 = \iota$.
The arguments above may be repeated verbatim to show that $\diam (H_t) \geq 1/4$ for all $t\in I$.
Indeed, $H$ may be lifted to a homotopy $\tilde{H}:\sph^2 \times I \to \Eqp$ such that $\tilde{H}_0 = i_0$.
Let $q\in \RP^2$, and $\gamma$ a non-trivial loop based at $q$ as before.
Again, $\gamma$ lifts to a path in $\sph^2$ connecting $x$ to $-x$, where $q$ is the class of $\{x, -x\}$.
Note that the following diagram commutes
\begin{equation*}
\begin{tikzcd}
\sph^2 \times I \arrow[r, "\tilde{H}"] \arrow[d, "\operatorname{prj}\times \mathds{1}"'] & \mathcal{E}q_+ \arrow[d, "f"] \\
\RP^2\times I \arrow[r, "H"]     & \mathcal{E}q                 
\end{tikzcd}
\end{equation*}
where, $\operatorname{prj}:\sph^2\to\RP^2$ is the covering map, and $f:\Eqp\to\Eq$ is the orientation forgetful map.
Hence, $\tilde{H}_t(x)$ and $\tilde{H}_t(-x)$, are the same unoriented equator, for all $t\in I$.
Since $i_0(x)$ and $i_0(-x)$ have opposite orientations, so does $\tilde{H}_t(x)$ and $\tilde{H}_t(-x)$ for all $t\in I$, by continuity.
Now, for any $t\in I$, the arguments before may be repeated to find a point $\tilde{y}\in \sph^2$, such that $d_H(\tilde{H}_t(\tilde{y}), \tilde{H}_t(x)) \geq 1/4$, and $d_H(\tilde{H}_t(\tilde{y}), \tilde{H}_t(-x)) \geq 1/4$.
Put $y=\operatorname{prj}(\tilde{y})$ the projection of $\tilde{y}$ to $\RP^2$.
We have thus found that,
\begin{equation*}
d_H(H_t(y), H_t(q)) \geq \frac{1}{4}.
\end{equation*}
Therefore, $\diam (H_t) \geq 1/4$ for all $t\in I$, as claimed.

The fourth and final elementary observation is that the naive bound of $1/2$ in \eqref{eq:naive-bound} is not optimal, as the following perturbation of $i_0$ has a strictly smaller diameter.
The idea of the construction is to perturb $i_0$ such that no equator in its image appears together with its oppositely oriented copy.
To bound the Hofer distance between perturbed equators we use the following lemma.
\begin{lemma}\label{lem:finite-intersection-lagrangians}
	Let $L,L'\subset\sph^2$ be two equators obtained as $C^2$-small perturbations of a great circle, taken with opposite orientations, such that $\# L \cap L' < \infty$.
	Then,
	\begin{equation*}
	d_H (L, L') < \frac{1}{2}.
	\end{equation*}
\end{lemma}

Before proving Lemma \ref{lem:finite-intersection-lagrangians}, let us first go into some of the details of the construction.
Observe that in a small neighborhood $U$ of a great circle $L$, one may apply a Hamiltonian diffeomorphism to perturb $L$ to a graph of a smooth function of the angle parameter of $L$.
Indeed, let $(q,p)$ be coordinates in $U$, such that $q$ is a parameter along $L$, and  $p$ a coordinate in an orthogonal direction.
Consider $\sph^1$ as $\sfrac{\RR}{2\pi\ZZ}$, and let $f(q):\sph^1\to(-\delta, \delta)$ be a continuous function with $f(0)=0$, and $\int_{\sph^1}f(t)dt=0$, where $\delta>0$ is such that the graph of $f$ is contained in $U$ when considered as a function along $L$.
Then, $H(q,p) = \int_0^q f(t)dt$ is a Hamiltonian, whose time 1 flow, $\phi_H^1$, satisfies that $\phi_H^1(L)$ is the graph of $f$.

Choose an open cover of $\sph^2$ by disjoint pairs of antipodal topological discs, denoted by $\ucal = \{U_1, V_1, \ldots, U_k, V_k\}$, where $U_j$ and $V_j$ are antipodal.
For ease of notation, put $L_p = i_0(p)$ for $p\in \sph^2$.
In each neighborhood $W\in \ucal$, pick a parametrization for all the great circles in $i_0(W)$.
Pick distinct prime numbers $r_j$, $s_j$ for each pair of $U_j, V_j\in \ucal$, and use the above parametrizations to define, for each $j=1,\ldots,k$, and each $p\in U_j$, $q\in V_j$, functions $f_p:L_p\to(-\delta, \delta)$, $g_q:L_q\to(-\delta, \delta)$, by $f_p(t) = \delta \sin(r_j t)$, $g_q(t) = \delta \sin(s_jt)$.
Using a partition of unity subordinate to $\ucal$, we merge these functions to $\hat{f}_p:L_p\to(-\delta,\delta)$, defined for all $p\in \sph^2$, such that it is $C^2$-small.
By construction, for all $p\in\sph^2$ it holds that the graphs of $\hat{f}_p$, and $\hat{f}_{-p}$ are $C^2$-small perturbation of the same unoriented great circle, with a finite number of intersection points.
Thus, Lemma \ref{lem:finite-intersection-lagrangians} applies.

By the second elementary observation above, the Hofer diameter of $i_0$ is realized by antipodal points.
The above perturbation results in an embedding $i:\sph^2\to\Eqp$, where each pair of points are at Hofer distance strictly less than $1/2$ apart, whence by compactness,
\begin{equation*}
\girth([i_0]) < 1/2,
\end{equation*}
as claimed.

\begin{proof}[Idea of the Proof of Lemma \ref{lem:finite-intersection-lagrangians}]
	
	Denote by $D_1$, $D_2$ the filling discs of $L$, $L'$ respectively, such that $D_1\cap D_2$ has area $\epsilon \ll 1/2$.
	Pick a point $p_0$ in $\sph^2$ outside of $D_1 \cup D_2$, and ``puncture" the sphere at $p_0$ to obtain the following picture in $\RR^2$; see Figure \ref{fig:perturbing-the-natural-embedding}.
	\begin{figure}[h]
		\centering
		\includegraphics[width=0.95\textwidth]{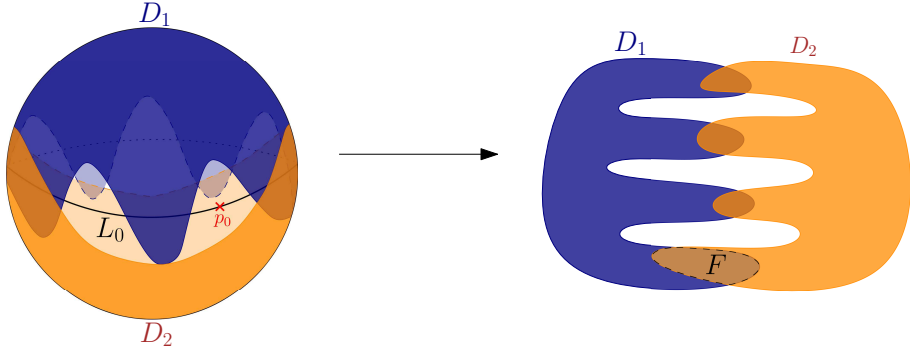}
		\caption{Perturbing equators to create an embedding with diameter smaller than $1/2$.}
		\label{fig:perturbing-the-natural-embedding}
	\end{figure}
	The two discs, $D_1$ and $D_2$, project to two regions in $\RR^2$, depicted in Figure \ref{fig:perturbing-the-natural-embedding} by the orange and blue blotches.
	By abuse of notation, denote these regions by $D_1$ and $D_2$ as well.
	Note that each of these has area $1/2$.
	By the hypothesis, the intersection $D_1\cap D_2$ has finitely many connected components; let $F$ be the one of largest area.
	Denote by $\epsilon'$ the area of $F$, so the area of $D_1\cap D_2 \setminus F$ is $\delta = \epsilon - \epsilon'$.
	Using techniques as in \cite{MR3395269}, which we demonstrate more thoroughly in Section \ref{sec:hofer-diam-calculation}, flow the smaller connected components into $F$, through $D_2$, such that their area is ``pushed" into $D_1$.
	A sketch of the flow is demonstrated in the Figure \ref{fig:dissolving-intersection-cnctd-cmps}, where there are only two connected components for the intersection.
	\begin{figure}[h]
		\centering
		\includegraphics[width=0.95\textwidth]{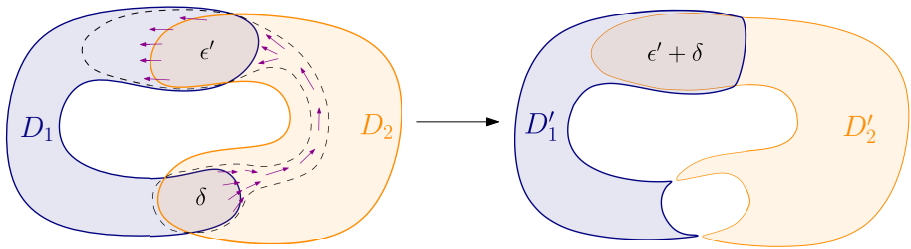}
		\caption{Flowing the intersection's connected components into the largest component through a flow described intuitively here by arrows.}
		\label{fig:dissolving-intersection-cnctd-cmps}
	\end{figure}
	Of course, this is possible since $D_1$ has area $1/2$, whereas the area flowed into $D_1$ is (slightly larger than) $\delta$, which is by assumption much smaller than $1/2$.
	Remark that an ``excess" area is flowed for technical reasons which we omit here.
	Since $D_1\cap D_2$ has a finite number of components, we may apply this procedure repeatedly until a single connected component remains.
	
	Denote by $D_1'$, and $D_2'$, the images of $D_1$, and $D_2$, respectively, under the procedure above.
	Observe the following few properties.
	First, the intersection $F':=D_1'\cap D_2'$ has a single connected component of area $\epsilon$. 
	Second, both $D_1'\setminus F'$ and $D_2'\setminus F'$ have area $1/2 - \epsilon$.
	Third and last, since $L$, $L'$, are graphs with respect to the same great circle, the resulting $D_1'$ and $D_2'$ each have only one connected component.
	Hence, after a change of coordinates we may consider $D_1'$, and $D_2'$, as two discs with the above areas and an intersection of area $\epsilon$.
	Since $\Area(D_1'\setminus F') = \Area(D_2'\setminus F')$, it follows that,
	\begin{equation*}
	d_H (L, L') \leq d_H(\partial D_1', \partial D_2') + \delta = \frac{1}{2} - \epsilon + \delta
	= \frac{1}{2} - \epsilon' < \frac{1}{2},
	\end{equation*}
	as claimed.
	
\end{proof}

\subsection{First look at pipe equators}\label{sec:embedding-construction}
\begin{figure}[h]
	\centering
	\includegraphics[width=0.4\textwidth]{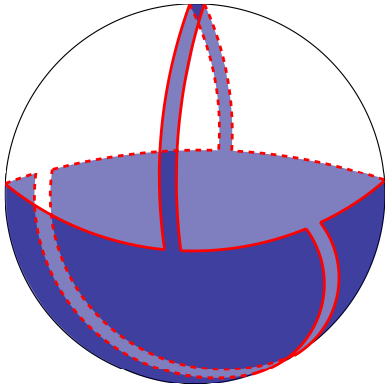}
	\caption{Pipe equator in red; the blue region denotes a connected component of the sphere with the equator removed.}
	\label{fig:piped-embedding}
\end{figure}

In this subsection we would like to introduce in rather loose terms the concept of ``pipe equators".
Let us begin with considering a cartoon of a typical equator we wish to examine, in Figure \ref{fig:piped-embedding}, and explain the motivation for the concrete construction which will follow in Sections \ref{sec:homotopy-construction} and \ref{sec:hofer-diam-calculation}.
As somewhat of a static picture, the red curve in Figure \ref{fig:piped-embedding} may be obtained by gluing two equal area slit discs to one another to form a sphere.
Choosing the discs to have area $1/2$ each, and the strips to have area $\delta$ each, results in a curve on the glued sphere, which separates it to two equal area connected components.
Therefore, this curve is in fact an orientable equator.
For the rest of the paper, these cut out strips will be called ``pipes", and equators formed this way will be called ``pipe equators".

Observe that these pipe equators are determined by three parameters: the area $\delta$ of the pipes (which we think of as fixed); the area $a$ to one fixed side of the upper pipe; and $b$ the area to one fixed side of the lower pipe. 
Remark that we chose to begin with this static picture in order to more easily explain the obstruction to the parameters determining pipe equators, which in turn leads to the homotopy constructed in Section \ref{sec:homotopy-construction}.
In practice, a more useful approach to these equators is dynamical in essence, and is investigated further in Section \ref{sec:hofer-diam-calculation}.

Let us briefly discuss these parameters and consider an unsuccessful yet instructive first attempt at constructing an embedding $\sph^2\to\Eqp$ via pipe equators.
Observe that a priori $a\in(0,1/2)$, since for $a=0,1/2$ the pipe ``collapses" on $L_0$ and so a pipe equator is undefinable.
Similarly, $b\in (0,1/2)$.
So, denote $I=(0,1/2)$ and consider $I^2$ as the pipe equators' parameters space. 
The above suggests the existence of an embedding $\iota:I^2\to \Eqp$, such that the boundary of $I^2$ should be mapped to $L_0$. 
Thus, one might naively attempt to collapse $\partial I^2$ to a point in order to obtain a map $\sph^2\to \Eqp$; however, this map remains ill defined.
Since the area on the other side of the pipe in the upper hemisphere is $1/2-a-\delta$, and like $a$ it is bounded in $(0,1/2)$, it follows that $a$ cannot be more that $1/2-\delta$; a similar argument shows that $b$ should lie in $(0, 1/2-\delta)$.

To properly define such a continuous embedding, we homotope between $i_0$ and an embedding $i_p:\sph^2 \to \Eqp$ with the following property; namely, it is obtained as a continuous transition from a constant map near the boundary $\partial I^2$, to $\iota$ away from it; see Section \ref{sec:homotopy-construction} for details about this construction.
Of course, $i_p$ belongs in the same homotopy class of $i_0$ in $\pi_2(\Eqp)$.
We show below that its Hofer diameter is bounded above by $1/3 + \delta$.
Essentially, this is attained by equators represented by points far from the boundary $\partial I^2$, i.e., by pipe equators.

\begin{structure}
	Section \ref{sec:homotopy-construction} constructs a homotopy between the natural embedding $i_0$ and an embedding whose image mostly consists of pipe equators.
	Then, Section \ref{sec:hofer-diam-calculation} shows that such an embedding has a Hofer diameter of $1/3$, up to an arbitrarily small parameter. 
	For completeness, Section \ref{sec:algebraic-topology} shows that $\pi_2(\Eqp)\neq 0$ and that the homotopy class of $i_0$ is non-trivial.
\end{structure}

\begin{acknowledgements*}
	I would like to thank Dr. M. Khanevsky for his contribution to this project, as well as for his patience and guidance.
\end{acknowledgements*}

\section{Homotopy between embeddings}\label{sec:homotopy-construction}

In this section we construct a homotopy between $i_0$, and an embedding $\sph^2\hookrightarrow\Eqp$ which depends on $0 < \delta \ll 1$ such that for most points on the sphere their image under the embedding is a pipe equator.
The construction goes through projecting to the planar unit disc, and applying a homotopy there.
Before going into the construction details, first fix some $0 < \epsilon, \delta \ll 1$, where $\delta$ denotes the area of a pipe in a pipe equator as in Subsection \ref{sec:embedding-construction}, and $\epsilon$ will be used below to demarcate a region of the equators' parameter space where the homotopy transitions from the constant map to a pipe equator embedding.
Let us now introduce a suitable parametrization of the sphere. 
It will be used throughout this section, starting with defining the projection to the unit disc.
Remark that in order to describe this parametrization slightly more easily, we reverse in this section the orientations of equators in the image of $i_0$.
This reversal does not interfere with the generality of our arguments.

\begin{figure}[H]
	\centering
	\includegraphics[width=0.8\textwidth]{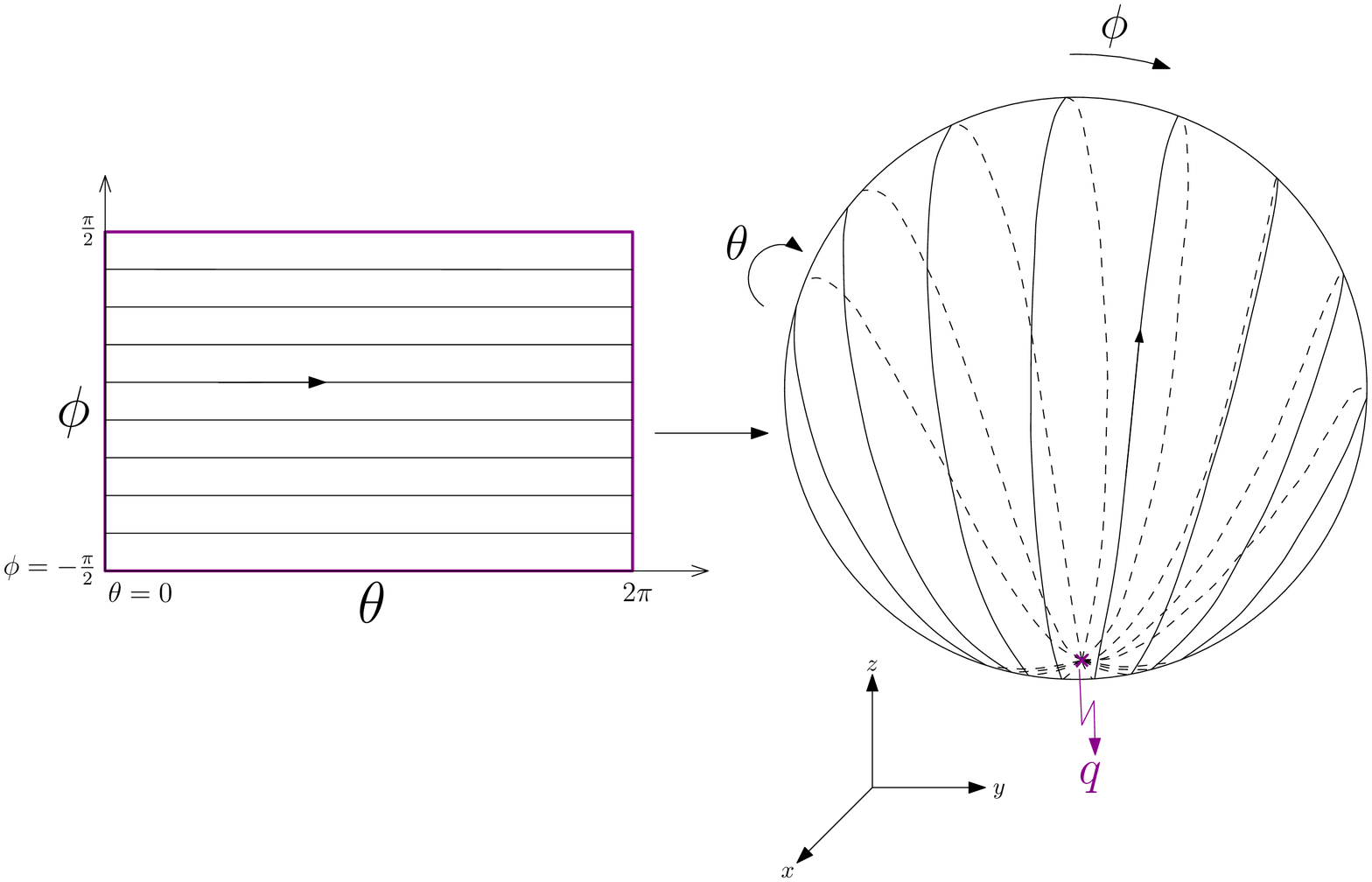}
	\caption{Parametrization of the sphere by angles.}
	\label{fig:sphere-fan-slice}
\end{figure}

Denote by $q\in\sph^2$ the south pole of the sphere.
Observe that every point on $\sph^2\setminus\{q\}$ lies on a circle through $q$ whose angle with respect to the $xz$-plane is $\phi\in(-\pi/2, \pi/2)$, and the point's location on the circle is given by $\theta\in(0, 2\pi)$, such that $\theta=0,2\pi$ refers to $q$; see Figure \ref{fig:sphere-fan-slice}.
Denote these circles by $C_\phi$; note that they are obtained as the intersection of $\sph^2$ with affine planes $\Pi_\phi$ forming an angle $\phi$ with the $xz$-plane and parallel to the $x$-axis.

Let us construct a projection to the planar unit disc.
Fix some $0\leq \phi<\pi/2 - \epsilon$, and consider the equators given by $i_0(x)$, where $x\in C_\phi$. 
These are rotations of $L_0$ about an axis, $l_\phi$, perpendicular to $\Pi_\phi$ such that it passes through the center of $\sph^2$, and its positive direction is determined by the intersection with the upper hemisphere. 
Indeed, each $x\in C_\phi$ determines a counterclockwise rotation $R=R(\theta, l_\phi)$ of $\sph^2$, with angle $\theta$ about the axis $l_\phi$, where $\theta=\theta(x)$ is the $\theta$ coordinate of $x$. 
Observe that as $L_0$ rotates, it forms an open cap in the upper hemisphere, $\dcal$, consisting of points which are unswept by the rotated $L_0$.
Note that $\dcal$ is the spherical cap formed about the axis $l_\phi$ in the positive direction, with an angle $\phi$ measured at the center of $\sph^2$, relative to $l_\phi$.
Denote the center of $\dcal$ by $x_\phi$.
As a counterpart to $\dcal$, its reflection with respect to $\Pi_\phi$ is a cap $\overline{\dcal}$ contained in the lower hemisphere, $\hcal_0$, and tangent to $L_0$, so it consists of points that are always contained in $\hcal_0$. 
Indeed, as $x$ ranges over $C_\phi$, the lower hemisphere, $\hcal_0$, and $\overline{\dcal}$, rotate via $R(\theta, l_\phi)$, for $\theta=\theta(x)$. 
Clearly, $R(\theta, l_\phi)( \overline{\dcal})$ is contained in $\hcal_\theta = R(\theta, l_\phi) (\hcal_0)$ and is tangent to $L_\theta = R(\theta, l_\phi) (L_0)$ for all $\theta\in(0,2\pi)$.

\begin{figure}[H]
	\centering
	\includegraphics[width=0.95\textwidth]{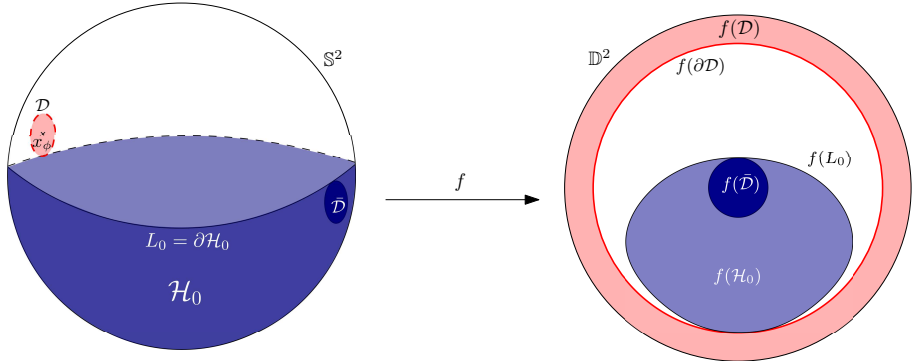}
	\caption{Projection of the sphere to the unit disc with some distinguished regions.}
	\label{fig:projection-to-disc-half-sphere}
\end{figure}

Puncture the sphere at $x_\phi$ to obtain a projection to the unit disc, i.e., a diffeomorphism $f:\sph^2\setminus\{x_\phi\}\to\DD^2$ to the open disc of area $1$, such that it preserves areas and translates rotations about $l_\phi$ to rotations about the center of the disc.
Note that the inverse of this diffeomorphism can be described by taking the quotient of the open unit disc by its boundary, thus collapsing it to a point.
The picture in the disc as appears in Figure \ref{fig:projection-to-disc-half-sphere}.
Remark that the rotation translation requirement on $f$ has the following interpretation in Figure \ref{fig:projection-to-disc-half-sphere}. 
As $x$ ranges over $C_\phi$ the rotation it induces on $\hcal_0$ translates under $f$ to a rotation of the light blue region in the disc such that it always contains the dark disc, tangent to it and to the red circle.

\begin{figure}[H]
	\centering
	\includegraphics[width=0.9\textwidth]{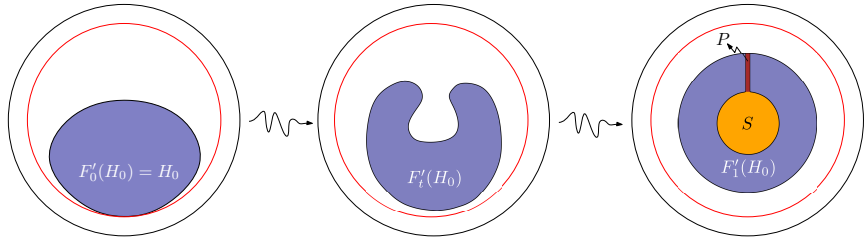}
	\caption{Homotopy of the disc.}
	\label{fig:homotopy-in-the-disc}
\end{figure}

Denote by $H_0\subset \DD^2$ the region $f(\hcal_0)$, it has half the area of the disc. 
Consider now the following geometric homotopy, associated with $\theta=0,2\pi$, defined as in Figure \ref{fig:homotopy-in-the-disc}, and denoted by $F'_t$.
That is, the homotopy is defined by deforming $H_0$, starting at the picture on the left to that on the right, such that it preserves areas; particularly, $\Area(F'_t(H_0))$ is independent of $t$. On the left $H_0$ remains unchanged as $F'_0$ is the identity map. On the right we have $F'_1(H_0)$ contained inside a slit annulus. Denote the slit by $P$ and the inner disc of the annulus by $S$, as in the figure.
The homotopy $F'$ may be constructed such that $\Area(P)=\delta$, and
\begin{equation*}
	\Area(S) = \frac{1-2\delta}{4\epsilon - 2\pi}\phi - \frac{\delta}{2} + \frac{1}{4}.
\end{equation*}
Particularly, $\Area(S)=0$ for $\phi = \pi/2 - \epsilon$, and $\Area(S) = 1/4-\delta/2$ for $\phi=0$.

\begin{figure}[H]
	\centering
	\includegraphics[width=0.7\textwidth]{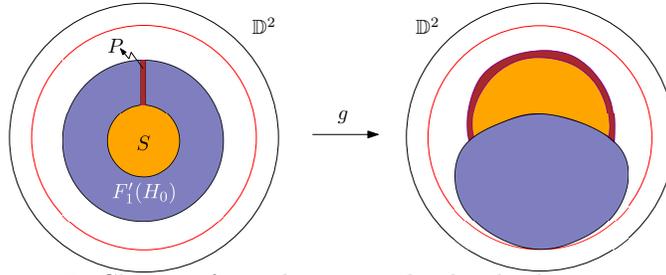}
	\caption{Change of coordinates in the disc leading to a projection of a pipe connected to a hemisphere.}
	\label{fig:change-of-coordinates-in-the-disc}
\end{figure}

\begin{figure}[H]
	\centering
	\includegraphics[width=0.35\textwidth]{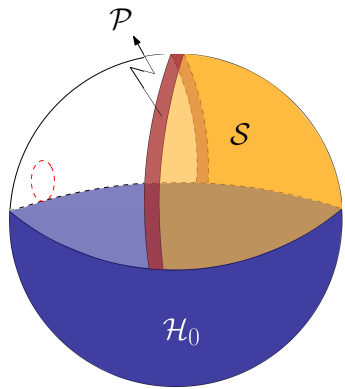}
	\caption{Hemisphere with a pipe attached to it.}
	\label{fig:sphere-initial-pipe-position}
\end{figure}

To define $F'_t$ for $\theta\in(0,2\pi)$, put
\begin{equation*}
	F'_{t,\theta} = R_{\theta} \circ F'_t \circ R_{-\theta},\quad t\in[0,1], \theta\in[0,2\pi],
\end{equation*}
where $R_{\alpha}$ is a positive rotation of the disc in an angle $\alpha$ with respect to the $x$-axis in the usual planar coordinates.
The upshot is that when $\theta$ goes over $(0,2\pi)$ the slit $P$ is rotated within the annulus. 
To make use of this behavior, we would like to apply a change of coordinates which translates the rotation in the annulus to a flow of the lower hemisphere through a pipe.
Namely, define the change of coordinates $g:\DD^2 \to \DD^2$ as in Figure \ref{fig:change-of-coordinates-in-the-disc}, such that $g(S)$ coincides with $f(\scal)$, where $\scal\subset \sph^2$ is the area bounded by the pipe, $\pcal=f^{-1}(P)$, in the upper hemisphere, as appears in Figure \ref{fig:sphere-initial-pipe-position}.
We further require $g$ to preserve areas.

\begin{figure}[H]
	\centering
	\includegraphics[width=0.7\textwidth]{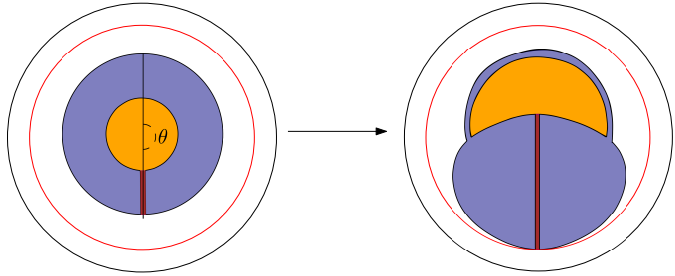}
	\caption{Rotation followed by a change of coordinates in the disc.}
	\label{fig:rotated-change-of-coords-in-disc}
\end{figure}

\begin{figure}[H]
	\centering
	\includegraphics[width=0.4\textwidth]{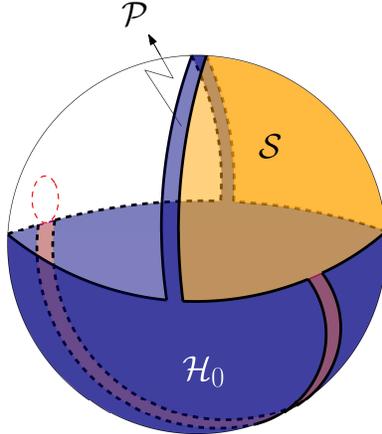}
	\caption{Pipe equator obtained by a projection to the unit disc and applying a homotopy there.}
	\label{fig:sphere-rotated-pipe-position}
\end{figure}

Then, as $\theta$ varies, we obtain the following picture after applying $g$; see Figure \ref{fig:rotated-change-of-coords-in-disc}.
By the well known Alexander trick, it holds that $g$ is homotopic to the identity, since it is compactly supported. 
Let $g_t$ be a homotopy, such that $g_0$ is the identity, $g_1=g$, and for all $t\in[0,1]$, it holds that $g_t$ preserves areas.
Note that $g_t$ can be chosen such that $\Area(g_t(S))$ is constant with respect to $t$.
Define
\begin{equation*}
F_{t,\theta} = g_t \circ F'_{t,\theta},
\end{equation*} 
Note that $F_{t,\theta}$ preserves the area of $H_0$ as $t$ and $\theta$ vary, so the boundary $\partial [f^{-1} \left(F_{t,\theta}(H_0)\right)]$ is an equator of $\sph^2$.
The equator obtained by this procedure for $t=1$ appears in Figure \ref{fig:sphere-rotated-pipe-position} as the bold black outline of the blue region (the dashed sections are in the back of the sphere). 
This equator is precisely a pipe equator, as constructed in Subsection \ref{sec:embedding-construction}. 

Now that we have all the pieces in place, we can construct the homotopy of embeddings, one hemisphere at a time.
Denote by $\sph^2_+:=\{(\theta, \phi) \in \sph^2 : \phi > 0\}$, and by $\Df$ the space of smooth simple closed loops in $\DD^2$ such that they divide the area of $\DD^2$ in half.
To define a homotopy $G:[0,1]\times \sph^2_+ \to \Eqp$ we require the following map $\tilde{f}:\Eqp \to \Df$ induced by $f$, namely $L\mapsto f(L)$. 

Define,
\begin{equation*}
G_t(\theta,\phi)  = \tilde{f}^{-1} \circ 
\begin{cases}
	F_{0,\theta} \circ \tilde{f} \circ i_0 (\theta,\phi),& \phi \in [0, \epsilon], \\
	F_{\left[t\cdot \frac{\phi - \epsilon}{\epsilon}, \theta\right]} \circ \tilde{f} \circ i_0 (\theta,\phi),& \phi \in [\epsilon, 2\epsilon], \\
	F_{t,\theta} \circ \tilde{f} \circ i_0 (\theta,\phi), & \phi \in [2\epsilon, \frac{\pi}{2} - 2\epsilon],\\
	F_{\left[t\cdot \frac{(\pi/2 - \epsilon) - \phi}{\epsilon}, \theta\right]} \circ \tilde{f} \circ i_0 (\theta, \phi),& \phi \in [\frac{\pi}{2} - 2\epsilon, \frac{\pi}{2} - \epsilon], \\
	F_{0,\theta} \circ \tilde{f} \circ i_0 (\theta,\phi),& \phi\in[\frac{\pi}{2}-\epsilon, \frac{\pi}{2}),
\end{cases}
\end{equation*}
where, by $F_{\bullet,\theta} \circ \tilde{f}\circ i_0$ we mean to apply $F_{\bullet,\theta}$ to the simple closed curve given by $\tilde{f} \circ i_0$. As discussed earlier, $F_{\bullet,\theta}$ preserves the area enclosed by $\tilde{f} \circ i_0$. Thus, $G_t$ is well defined. 
Evidently, $G_0 \equiv i_0$, and $G_1$ embeds $\sph^2_+$ in $\Eqp$, up to $\epsilon$ and $\delta$.
Note that $\epsilon$ is chosen such that $G$ transitions continuously from $i_0$ to the pipe equator embedding. 
Such a transition is required since the latter embedding cannot be defined for sufficiently small $\phi$, due to the area required for the pipe, as explained in Subsection \ref{sec:embedding-construction}.

For $\phi \in (-\pi/2, 0)$, the above geometric picture is mirrored, and so the construction may be mirrored as well. As $\phi \searrow 0$, both $\dcal$ and $\overline{\dcal}$ shrink in area. For $\phi=0$ both caps collapse to a point. Then as $\phi$ descends below $0$, $\dcal$ becomes the cap which is always covered by hemispheres as $\theta$ varies, whereas $\overline{\dcal}$ becomes the cap which is never obtained by the hemispheres.
Now, pick $x_\phi$ to be the center of $\overline{\dcal}$, puncture the sphere at $x_\phi$ to obtain $h:\sph^2\setminus\{x_\phi\}\to\DD^2$, and $\tilde{h}:\Eqp\to \Df$, analogously to $f$ and $\tilde{f}$ above.

Define $H:[0,1]\times \sph^2_{-} \to \Eqp$, by
\begin{equation*}
H_t(\theta,\phi)  = \tilde{h}^{-1} \circ 
\begin{cases}
	F_{0,\theta} \circ \tilde{h} \circ i_0 (\theta,\phi),& \phi \in [-\epsilon, 0], \\
	F_{\left[-t\cdot \frac{\epsilon + \phi}{\epsilon}, \theta\right]} \circ \tilde{h} \circ i_0 (\theta,\phi),& \phi \in [-2\epsilon, -\epsilon], \\
		F_{t,\theta} \circ \tilde{h} \circ i_0 (\theta,\phi),& \phi \in [-\frac{\pi}{2} + \epsilon, -2\epsilon],\\
	F_{\left[t\cdot \frac{\phi + (\pi/2 - \epsilon)}{\epsilon}, \theta\right]} \circ \tilde{h} \circ i_0 (\theta, \phi),& \phi \in [-\frac{\pi}{2} + 2\epsilon, -\frac{\pi}{2} + \epsilon], \\
	F_{0,\theta} \circ \tilde{h} \circ i_0 (\theta,\phi),&\phi\in(-\frac{\pi}{2}, -\frac{\pi}{2} + 2\epsilon],
\end{cases}
\end{equation*}
where $\sph^2_{-} := \{(\theta, \phi)\in \sph^2 : \phi < 0\}$.
Now we may define the following amalgamation of $G$ and $H$ as $A:[0,1]\times \sph^2 \to \Eqp$, by
\begin{equation*}
A_t (\theta, \phi) = \begin{cases}
H_t(\theta, \phi),& \phi \in (-\frac{\pi}{2}, 0], \\
G_t(\theta, \phi),& \phi \in (0, \frac{\pi}{2}).
\end{cases}
\end{equation*}
Observe that for all $t\in[0,1]$ and all $\phi \in (-\epsilon, \epsilon)$ it holds that $A_t \equiv i_0$. Since $i_0$ is continuous, so is $A$ for these values. For all other $\phi$, $A$ coincides with either $G$ or $H$ and is therefore continuous.

This completes the construction of the homotopy between the natural embedding $i_0$ and a pipe equator embedding $i_p$ associated with $\delta$ and $\epsilon$.

\section{The Hofer diameter of a pipe equator embedding}\label{sec:hofer-diam-calculation}
In this section we provide an upper bound on the Hofer distance between two pipe equators, represented by points away from the boundary of the parameters space.
This is done by considering a slightly more convenient planar setting, where two types of flows are applied.
The upper bound is then given by calculating the energy cost of applying these flows in two possible ways.

First, let us describe the dynamical approach to the formation of pipe equators, which will then be used as part of the planar picture.
Consider $\sph^2$, and construct a thin strip of area $0<\delta \ll 1/2$ in the upper hemisphere, such that the strip connects $L_0$ to itself.
As before, we refer to this strip as a pipe, and denote the area to one of its sides in the upper hemisphere by $a$.
Note that while the choice of a side whose area we denote is arbitrary, it does not affect the following, as $a$ is uniquely determined by the area on the other side of the pipe.

\begin{figure}[h]
	\centering
	\includegraphics[width=0.45\textwidth]{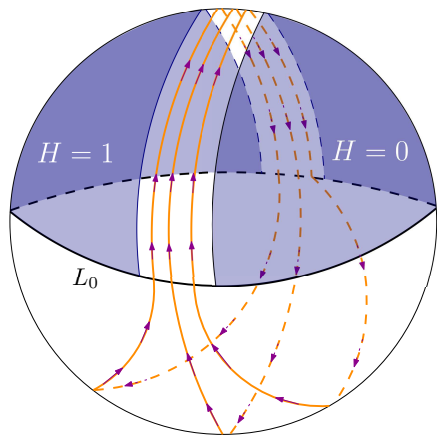}
	\caption{Hamiltonian diffeomorphism creating a pipe equator, and a few flow lines.}
	\label{fig:placeholder-Ham-flow-on-the-sphere}
\end{figure}

Given such a pipe, we flow the lower hemisphere through it via a Hamiltonian diffeomorphism, reasoning similarly to constructions given in Sections 2.1 and 2.2 in \cite{MR3395269}.
Let $H$ be an autonomous Hamiltonian such that on one of the pipe's sides $H=0$, on its other side $H=1$, and $H$ is linearly interpolated in between; see Figure \ref{fig:placeholder-Ham-flow-on-the-sphere}, where the regions in blue are where $H$ is constant.
Apply a $C^\infty$ smoothing to $H$ near the singular points. 
Denote by $\phi^t_H$ the time-$t$ map of the flow of $H$. 
As explained in \cite{MR3395269}, after time $t+\epsilon$, a region of the lower hemisphere of area $t$ will be flowed along the pipe, where $\epsilon$ depends on the choice of smoothing of $H$ and can be made arbitrarily small.
Given that the lower hemisphere has been flowed for time $t+\epsilon>\delta$, the area of the pipe is pushed in the lower hemisphere to form a strip of area $\delta$, which in our terminology is also called a pipe.
Observe that the time parameter of the flow uniquely determines the area of the lower hemisphere not yet flowed through the pipe by $b=1/2 - t -\delta$. 
This area is the second parameter we mentioned of the resulting pipe equator. 
Hence, in the dynamical approach we define the pipe equator associated with $a,b$, as the image of $\phi^t_H(L_0)$, where $t = t(b)$, and $H$ depends on $a$.

Note that in Section \ref{sec:homotopy-construction} we used angle parameters to describe pipe equators, whereas here we use area parameters.
By the construction of the pipe equator embedding, we may translate one to the other.
Observe that $a=\Area(S)$ in the notation of Section \ref{sec:homotopy-construction}, for which we provided a formula for $\phi \in (0, \pi/2-\epsilon)$.
Hence, combining with the corresponding formula for negative $\phi$, we obtain,
\begin{equation*}
	a(\phi) = \frac{2\delta - 1}{2\pi - 4 \epsilon} \phi + \frac{1}{4} - \frac{\delta}{2},
\end{equation*}
As mentioned, $b$ is the area not yet flowed through the pipe, and by construction of the pipe equator embedding is chosen to depend on $\theta$ linearly; namely,
\begin{equation*}
	b(\theta) = \frac{2\delta - 1}{4\pi}\theta + \frac{1}{2} - \delta.
\end{equation*}

\begin{figure}[H]
	\centering
	\includegraphics[width=0.5\textwidth]{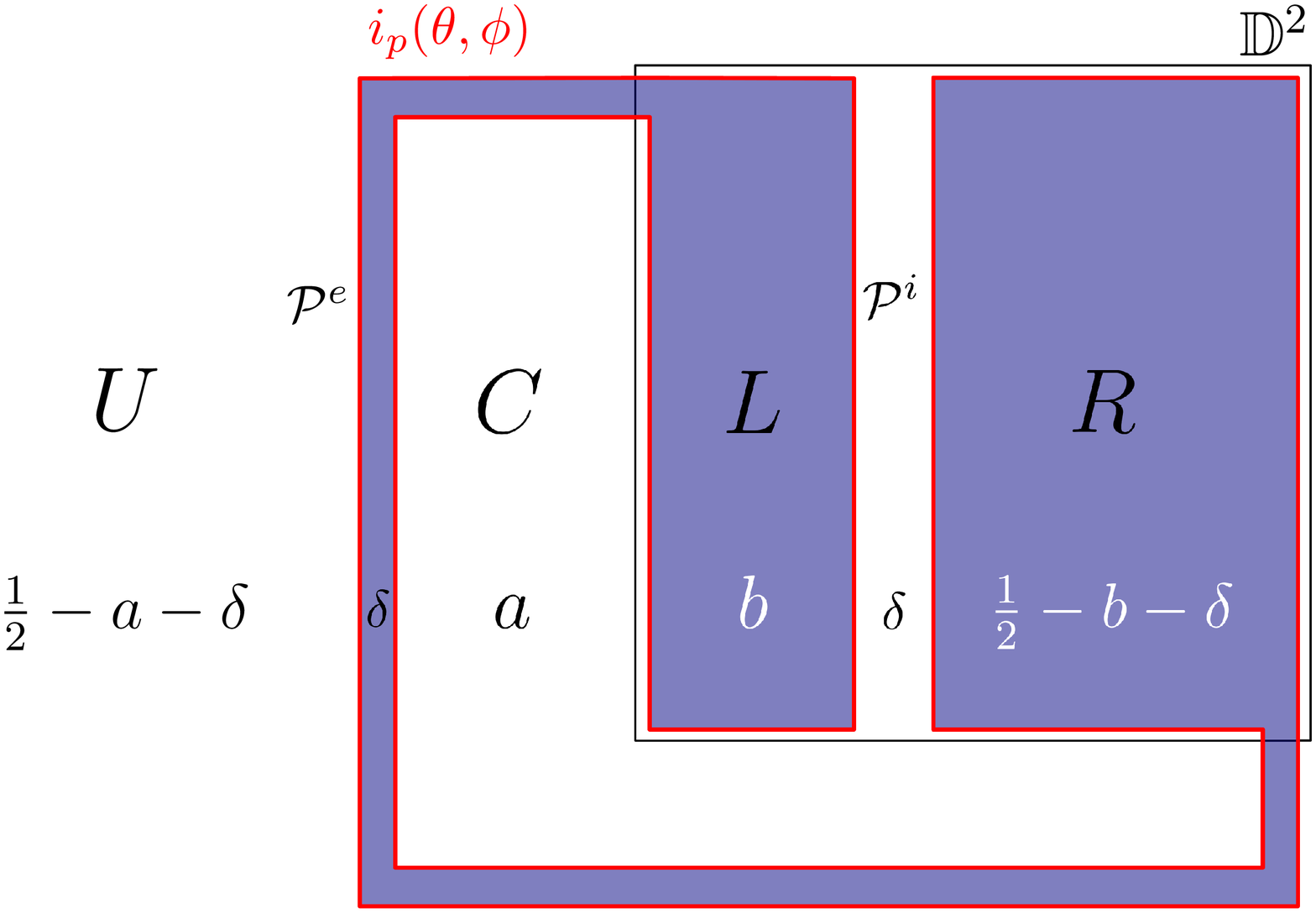}
	\caption{Pipe equator projected to the plane. Note the colors match Figure \ref{fig:piped-embedding}.}
	\label{fig:pipes-in-cube-initial-position}
\end{figure}

Let us now describe the planar outlook.
In the notation of the previous section, since we fixed $\delta$ and $\epsilon$, the parameter $\phi$ is bounded away from $0$ and there is some cap, $\dcal$, unswept by the flow of $L_0$.
Observe that $\ucal := \sph^2 \setminus \overline{\dcal}$ is open, and may be considered as a neighborhood containing all pipe equators associated with $\delta$.
By restricting the discussion to $\ucal$ and changing coordinates, we may consider pipe equators as curves in the plane, as in Figure \ref{fig:pipes-in-cube-initial-position}.
For ease of notation we denote this restriction to $\ucal$ and change of coordinates as a map $q:\ucal \to \RR^2$.
Let us clarify the notations of Figure \ref{fig:pipes-in-cube-initial-position}.
Let $\gamma$ be a pipe equator as the image of $L_0$ under a Hamiltonian flow $\phi$, such that its image under $q$ is the red curve in the figure. 
Denote by $\hcal_0$ the lower hemisphere of $\sph^2$, and $\pcal$ the pipe in the upper hemisphere.
Denote by $\hcal_0^r$ the region of $\hcal_0$ that has not already been flowed through $\pcal$, and by $\hcal_0^f$ the region that has been flowed.
Thus, in terms of Figure \ref{fig:pipes-in-cube-initial-position}, $L$ denotes the image $q(\hcal_0^r)$, similarly $R$ denotes $q(\hcal_0^f)$.
Likewise, denote the upper hemisphere of $\sph^2$ by $\hcal_1$, and analogously $\hcal_1^f$, and $\hcal_1^r$, as the regions that have (and have not) been flowed through the pipe in the lower hemisphere.
Apply the change of coordinates $q$, so $C$ denotes $q(\hcal_1^r)$, and $U$ denotes $q(\hcal_1^f)$ in terms of Figure \ref{fig:pipes-in-cube-initial-position}.
Finally, let $\pcal^e$ denote the image under $q$ of the pipe in the upper hemisphere, and $\pcal^i$ denote the image under $q$ of the pipe in the lower hemisphere.
As the pipe equator $\gamma$ is determined by the areas $a$ and $b$, so are the areas in the plane. 
Namely, $L$ has area $b$, $C$ has area $a$, the pipes both have area $\delta$, and $R$ and $U$ have complementing areas to $L$, $C$, and the pipes.

\begin{figure}[H]
	\centering
	\includegraphics[width=0.45\textwidth]{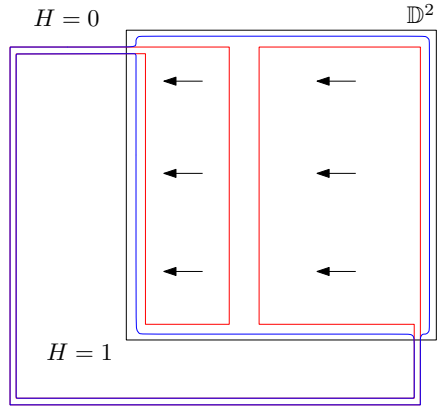}
	\caption{Flowing equators in the disc.}
	\label{fig:flowing-pipes-in-disc}
\end{figure}

\begin{figure}[H]
	\centering
	\includegraphics[width=0.45\textwidth]{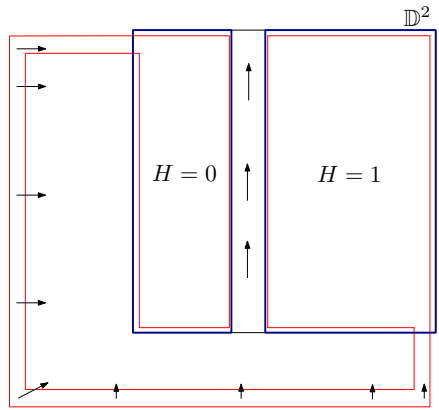}
	\caption{Shrinking the areas enclosed by the pipe and the disc.}
	\label{fig:flowing-areas-of-pipes}
\end{figure}

Therefore, given two pipe equators $\gamma_1=i_p(a_1, b_1)$, $\gamma_2 = i_p(a_2, b_2)$ as above, consider their Hofer distance
\begin{equation}\label{eq:Hofer-d-two-pipe-eqs-gen-pos}
d_H(i_p(a_1,b_1), i_p(a_2,b_2)),
\end{equation}
as a function of $0 < a_j, b_j < 1/2$, for $j=1,2$.
In the following we bound this distance from above by estimating the symplectic energy it would take to deform one to the other under $q$.
At the beginning of this section, we used arguments from \cite{MR3395269} to estimate the energy it takes to construct a pipe equator on the sphere.
In fact, the same arguments can be used to deform equators to one another by flowing the regions they enclose.
We consider two specific such deformations in the plane.
First, apply a flow whose Hamiltonian is as depicted in Figure \ref{fig:flowing-pipes-in-disc}, in order to flow the region $R$ towards $L$ in terms of Figure \ref{fig:pipes-in-cube-initial-position}. This flow may also be applied in negative time to flow $L$ towards $R$.
The second deformation we consider, consists in flowing the regions $C$ and $U$, in terms of Figure \ref{fig:pipes-in-cube-initial-position}.
By applying a flow associated with a Hamiltonian described in Figure \ref{fig:flowing-areas-of-pipes}, the region outside the disk $\DD^2$ is flowed through the pipe $\pcal^i$, such that the external pipe $\pcal^e$ is deformed to enclose a smaller area between $\DD^2$ and itself.
Intuitively speaking, this flow shrinks the area between the external pipe and the disc.

\begin{figure}[H]
	\centering
	\includegraphics[width=0.4\textwidth]{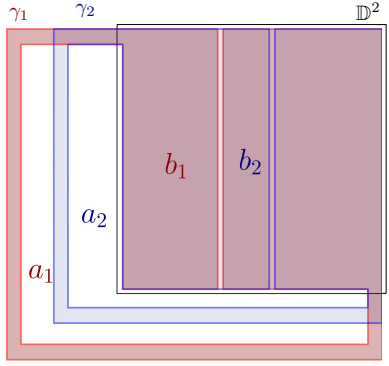}
	\caption{Two pipe equators in general position projected to the plane.}
	\label{fig:pipes-in-cube-two-general-positions}
\end{figure}

Given two equators as in Figure \ref{fig:pipes-in-cube-two-general-positions} we use the above deformations in order to describe two approaches to deform one equator into the other.
Roughly speaking, the first approach is to flow each of the two equators such that they both coincide with $\partial \DD^2$. The second approach is to flow the overlapping regions between the equators to make them coincide.

In the first approach, consider each pipe equator separately and observe that there are four different manners in which an equator can be flowed through the pipes such that it coincides with $\partial \DD^2$. In the notation of Figure \ref{fig:pipes-in-cube-initial-position}: the first (and second) possibility is to flow $L$ ($R$) through $\pcal^e$ and then flow $\pcal^e$ itself; the third (and fourth) possibility is to flow $C$ ($U$) through $\pcal^i$ and then flow $\pcal^i$ itself.
Since both pipes have area $\delta$, and in all four cases the pipes are also flowed, we calculate the required energy only up to $\delta$ and do not include the pipes in the arguments below.
Flowing $L$ is done with energy cost of $b$, as explained above; likewise, $R$ costs $1/2 - b$ to flow, $C$ costs $a$, and $U$ costs $1/2-a$.
Since we are interested in minimizing the cost, this yields the following upper bound on deforming a given pipe equator $i_p(a,b)$ to $\partial \DD^2$,
\begin{equation*}
	\min\left\{a, \frac{1}{2} - a, b, \frac{1}{2} - b \right\}.
\end{equation*}
Therefore, it follows that the first approach yields the following cost function,
\begin{equation*}
f_1(a_1, b_1, a_2, b_2) = \min\left\{a_1, \frac{1}{2} - a_1, b_1, \frac{1}{2} - b_1 \right\} + \min\left\{ a_2, \frac{1}{2} - a_2, b_2, \frac{1}{2} - b_2 \right\}.
\end{equation*}

In the second approach, instead of deforming both equators to coincide with $\partial \DD^2$, we deform them to coincide with each other.
This consists in flowing $L_1$ to coincide with $L_2$, or the other way around; and in flowing $C_1$ to coincide with $C_2$, or conversely $C_2$ to coincide with $C_1$.
Of course, flowing $L_1$ to $L_2$ is the same as $R_2$ to $R_1$, and the same holds for the relation between the $C$'s and $U$'s.
Hence, the cost of the former is $|b_1 - b_2|$, while the cost of the latter is $|a_1 - a_2|$.
Overall, the second approach yields the following cost function,
\begin{equation*}
f_2(a_1, b_1, a_2, b_2) = |a_1 - a_2| + |b_1 - b_2|.
\end{equation*}

So, altogether given two equators $L_1=i_p(a_1, b_1)$, $L_2(a_2, b_2)$, the following holds for their Hofer distance
\begin{equation*}
	d_H(L_1, L_2) \leq \min\{f_1(a_1, b_1, a_2, b_2), f_2(a_1, b_1, a_2, b_2)\}.
\end{equation*}
Therefore, put
\begin{equation*}
F := \min \{f_1, f_2\}.
\end{equation*}
It thus follows that
\begin{equation*}
	\diam(i_p) \leq \max F,
\end{equation*}
whence it remains to maximize $F$.

A maximum of $F$ in $[0,1/2]^4$ is an upper bound for \eqref{eq:Hofer-d-two-pipe-eqs-gen-pos}, up to $\delta>0$. We claim that $F \leq 1/3$, and that this value is attained at
\begin{equation*}
x_0=(1/3,1/3,1/6,1/6).
\end{equation*}
Clearly, $F(x_0) = 1/3$, whence $x_0$ is a maximal point for $F$.
Note that $x_0$ is not unique, there are at least $3$ other permutations of the coordinates of $x_0$ which yield a maximum for $F$.

To see that $F\leq 1/3$, first consider $x=(a_1,b_1,a_2,b_2)\in[1/6, 1/3]^4$.
It holds that $|a_1-a_2| \leq 1/6$, and likewise $|b_1 - b_2| \leq 1/6$, so $f_2(x) \leq 1/3$. Therefore,
\begin{equation*}
F(x) = \min\{f_1(x), f_2(x)\} \leq f_2(x) \leq \frac{1}{3}.
\end{equation*}
Second, for $x=(a_1,b_1,a_2,b_2)$ outside of $[1/6, 1/3]^4$, consider $f_1(x)$. It holds that 
\begin{equation*}
m_j(x) := \min \left\{a_j, \frac{1}{2} - a_j, b_j, \frac{1}{2} - b_j \right\} \leq \frac{1}{6},
\end{equation*}
for $j=1,2$. Indeed, for each $j=1,2$, pick any of the terms in the minimum expression above, and denote it by $c$. By our assumption, either $c\leq 1/6$ or $1/2 - c \leq 1/6$. Therefore, $m_j(x) \leq 1/6$, and so
\begin{equation*}
F(x) = \min\{ f_1(x), f_2(x) \} \leq f_1(x) = m_1(x) + m_2(x) \leq 1/3,
\end{equation*}
as claimed.

Finally, remark that all ``non-pipe" equators described by the pipe equator embedding (i.e., $i_p(\theta, \phi)$ for $\phi \in [-\pi/2, -\pi/2+\epsilon]\cup [\pi/2 - \epsilon, \pi/2]$) are, by construction, $\epsilon$-close to $L_0$.
Hence, up to $\epsilon$, the estimates above hold for them as well.
This finishes the estimate of the upper bound in Theorem \ref{thm:main-theorem}.

\section{Non-triviality of $[i_0]$}\label{sec:algebraic-topology}
Observe that $\Ham{\sph^2}$ acts on $\Eqp$ since Hamiltonian diffeomorphisms preserve areas.
Denote $S= \Stab_{\Ham{\sph^2}}(L_0)$ the stabilizer of the standard equator $L_0$ with positive orientation.
We have the following fibration
\begin{equation}\label{eq:fibration-of-Ham-over-Eqp}
\begin{tikzcd}
S \arrow[r] & \operatorname{Ham}(\mathbb{S}^2) \arrow[d] \\
& \mathcal{E}q_+                            
\end{tikzcd}
\end{equation}
This is true since $\Ham{\sph^2}$ is a Lie group, and so $\Ham{\sph^2} \to \Ham{\sph^2} / S$ is a fibration with fiber $S$; see, e.g. \cite{MR1688579} Sections 7.3-7.5. Here we used the fact that $\Eqp$ is homeomorphic to the orbit space $\Ham{\sph^2} / S$; this is true because the action is transitive, indeed every equator is obtained from $L_0$ by applying some Hamiltonian diffeomorphism.

Given the fibration above, we can consider the long exact sequence it induces on homotopy groups (see, e.g., Theorem 4.41 in \cite{MR1867354}),
\begin{multline}\label{eq:LES-Eqp}
	\pi_2(\Ham{\sph^2},\Id) \to  \pi_2(\Eqp,L_0) \xrightarrow{\partial} \pi_1(S,\Id) \\ \to \pi_1(\Ham{\sph^2},\Id) \to  \pi_1(\Eqp,L_0) .
\end{multline}

It holds that $\Ham{\sph^2}$ is homotopy equivalent to $SO(3)$. This is due to the fact that $\Ham{\sph^2} = \Symp_0(\sph^2)$ (see, e.g., Section 1.4 in \cite{MR1826128}), and since $\Diffp{\sph^2} \supset \Symp_0(\sph^2)$ has $SO(3)$ as a strong deformation retract, by Theorem A in \cite{MR112149}.
Since $SO(3)$ has the homotopy type of $\RP^3$, it follows that,
\begin{equation*}
	\pi_2(\Ham{\sph^2},\Id) \cong \pi_2(\Ham{\sph^2}) \cong \pi_2(\RP^3) = 0,
\end{equation*}
and similarly,
\begin{equation*}
	\pi_1(\Ham{\sph^2},\Id) \cong \pi_1(\RP^3) = \ZZ / \ZZ_{2}.
\end{equation*}

As for the $\pi_1(S,\Id)$ in \eqref{eq:LES-Eqp}, consider the following two maps. 
First, $i:\sph^1 \to S$ given by $t\mapsto R(2\pi t, \hat{z})$, where $R$ is described in the notation of Subsection \ref{sec:Elem-observations}, $\hat{z}$ is a unit vector in the direction of the $z$-axis, and $\sph^1=\RR / \ZZ$.
The second map is the evaluation $\eval:S \to \sph^1$ given by $\phi\mapsto l^{-1} (\phi(l(0)))$, where $l:\sph^1 \to L_0$ is some fixed parametrization of $L_0$. 
We have the following commutative diagram:
\begin{equation*}
\begin{tikzcd}
{\sph^1} \arrow[r, "i"] \arrow[rr, "\Id"', bend right] & {S} \arrow[r, "\eval"] & {\sph^1}
\end{tikzcd}
\end{equation*}
This diagram descends to a diagram of fundamental groups. 
\begin{equation*}
\begin{tikzcd}
{\pi_1(\sph^1,0)} \arrow[r, "i_*"] \arrow[rr, "\Id"', bend right] & {\pi_1(S,\Id)} \arrow[r, "\eval_*"] & {\pi_1(\sph^1,0)}
\end{tikzcd}
\end{equation*}
Hence it follows that $i_*$ is injective, and so $\pi_1(S, \Id)$ contains a $\ZZ$ subgroup.

Going back to \eqref{eq:LES-Eqp}, we have the following excerpt of an exact sequence:
\begin{equation}\label{eq:stab-exact-sequence}
\begin{tikzcd}
{\pi_2(\Eqp,L_0)} \arrow[r, "q"] & {\pi_1(S,\Id)} \arrow[r, "p"] & \ZZ/\ZZ_2
\end{tikzcd}
\end{equation}
Since $\ZZ < \pi_1(S,\Id)$, the map $p$ in \eqref{eq:stab-exact-sequence} cannot be injective. By exactness at $\pi_1(S,\Id)$ it follows that $\operatorname{Image}(q) \neq 0$.
This shows that $\pi_2(\Eqp, L_0) \neq 0$.

Finally, observe that $[i_0]\neq 0 \in \pi_2(\Eqp, L_0)$.
To see this, note that its image under the connecting morphism $\partial:\pi_2(\Eqp, L_0)\to\pi_1(S, \mathds{1})$ in \eqref{eq:LES-Eqp} is also non-trivial.
Indeed, recall that by \eqref{eq:fibration-of-Ham-over-Eqp} we have a fibration $\Ham{\sph^2}\to\Eqp$.
One may consider $i_0$ as a map $(\DD^2, \partial \DD^2)\to(\Eqp, L_0)$, by collapsing the $\sph^1$ boundary of $\DD^2$ to a point.
By the homotopy lifting property of fibrations, it follows that there exists a lift $\tilde{i}_0:(\DD^2, \partial \DD^2)\to (\Ham{\sph^2}, S)$.
Since the boundary $\partial \DD^2$ of the disk maps to $L_0$, it follows that the image of  $\tilde{i}_0|_{\partial \DD^2}$ is contained in the fiber $S$.
Therefore, one may consider $\tilde{i}_0|_{\partial \DD^2}$ as a map $(\sph^1, 0)\to (S, \Id)$.
By definition of the connecting morphism, it holds that 
\begin{equation}\label{eq:image-of-connecting-maps}
	\partial[i_0]= [\tilde{i}_0|_{\partial \DD^2}].
\end{equation}

We claim now that $[\tilde{i}_0|_{\partial \DD^2}]$ is not a unit, and by \eqref{eq:image-of-connecting-maps} it follows that $[i_0]$ is not a unit either.
Observe that a lift of $i_0$ to $\Ham{\sph^2}$ may be chosen to have its image contained in $\SO(3)$.
Since $\SO(3)$ is diffeomorphic to $S\sph^2$, the unit circle bundle in $T\sph^2$, it follows that a choice of a lift for $i_0:\DD^2 \to \SO(3)$ is equivalent to a choice of a section $v\in \Gamma(S\sph^2)$, such that it has a singular point at $N$, the north pole of $\sph^2$.
After multiplying $v$ by a function $\lambda:\sph^2\to\RR$ which vanishes exactly at $N$, we may apply the renowned Poincar\'e-Hopf Theorem, to obtain
\begin{equation*}
	\operatorname{index}_N(\lambda\cdot v) = \chi(\sph^2) = 2,
\end{equation*}
where $\chi(\sph^2)$ is the Euler characteristic of the sphere, which is well known to be 2.
Now, recall that $\eval_*$ denotes the map $\pi_1(S,\Id)\to \pi_1(\sph^1, 0)$ induced by evaluation.
By definition of the index, and considering $\pi_1(\sph^1, 0) \cong \ZZ$, we find that
\begin{equation*}
	\eval_* [\tilde{i}_0|_{\partial \DD^2}] = 2.
\end{equation*}
Notably, it follows that $[\tilde{i}_0|_{\partial \DD^2}] \neq 0$, as claimed.

\bibliography{refs}

\end{document}